\documentclass[submission,copyright]{eptcs}
\providecommand{\event}{QPL 2014} 
\usepackage{breakurl}             

\usepackage{graphicx}
\usepackage{amsmath}
\usepackage{alltt}
\usepackage{mathrsfs} 
\usepackage{amsfonts}
\usepackage{amssymb}%
\usepackage{amsthm}
\usepackage{amscd}
\usepackage{url}
\usepackage[all,knot,poly]{xypic} 
\usepackage{endnotes}
\usepackage{tikz}
\usetikzlibrary{plotmarks}
\usetikzlibrary{decorations.markings}

\usepackage{enumerate} 
\usepackage{comment} 

\newcommand{\mcC}{\mathcal{C}}

\newcommand{\scrT}{\mathscr{T}}

\newcommand{\C}{\mathbb{C}}
\newcommand{\R}{\mathbb{R}} 
\newcommand{\N}{\mathbb{N}}

\newcommand{\FF}{\mathbb{F}}

\newcommand{\x}{\mathbf{x}}

\newcommand{\size}{\operatorname{size}}

\newcommand{\ot}{\otimes}
\newcommand{\otc}{\otimes \cdots \otimes}

\newcommand{\ra}{\rightarrow}

\newcommand{\rae}{\!\rightarrow\!}

\newcommand{\nbhd}{\operatorname{nbhd}}

\newcommand{\id}{\operatorname{id}}

\newcommand{\Ob}{\operatorname{Ob}}

\newcommand{\Mor}{\operatorname{Mor}}

\newcommand{\cod}{\operatorname{cod}}
\newcommand{\dom}{\operatorname{dom}}

\newcommand{\Vect}{{\rm Vect}}

\newcommand{\Mat}{\operatorname{Mat}}
 
\newcommand{\BoolRel}{{\rm BoolRel}}

\newcommand{\FinRel}{{\rm FinRel}}
\newcommand{\NFinRel}{\N{\rm FinRel}}

\newcommand{\cpy}{\operatorname{copy}}

\newcommand{\coev}{\operatorname{coev}} 
\newcommand{\ev}{\operatorname{ev}}

\theoremstyle{plain}
\newtheorem{theorem}{Theorem}[section]

\newtheorem{proposition}[theorem]{Proposition}

\theoremstyle{definition}

\theoremstyle{definition}
\newtheorem{defn}[theorem]{Definition}
\theoremstyle{definition}

\theoremstyle{definition}
\newtheorem{example}[theorem]{Example}
\theoremstyle{definition}

\specialcomment{issue}
{\begingroup\color{red}\footnotesize}{\endgroup}

\excludecomment{baa}
\excludecomment{issue} 
\excludecomment{comment}

\title{Belief propagation in monoidal categories}
\author{Jason  Morton
\institute{Department of Mathematics\\Pennsylvania State University\\University Park, PA}
\email{morton@math.psu.edu}
}
\def\titlerunning{Belief propagation in monoidal categories}
\def\authorrunning{J. Morton}
\begin{document}
\maketitle

\begin{abstract}
We discuss a categorical version of the celebrated belief propagation algorithm.  This provides a way to prove that some  algorithms which are known or suspected to be analogous, are actually identical when formulated generically.  It also highlights the computational point of view in monoidal categories.
\end{abstract}

\section{Introduction}

We discuss a categorical version of the celebrated belief propagation algorithm \cite{pearl1982reverend}.  
This provides a way to prove that some  algorithms which are known or suspected to be analogous, are actually identical when formulated generically.  It also highlights the computational point of view in monoidal categories. The approach could also make possible software implementations that use a categorical formulation of algorithms to enable generic programming. 

The setting for this algorithm, and for many computational questions concerning monoidal categories, is a {\em diagram}  (Def. \ref{def_diagram}).   By a diagram we mean an  equivalence class of monoidal words over a finite tensor scheme, usually with certain additional properties.
An interpretation \cite{selinger2009survey} assigns values to each of the variables in the tensor scheme in a functorial way. 
Basic questions to 
ask of such a diagram, once interpreted in a particular category, include
\begin{enumerate}
\item compute a (possibly partial) contraction,
\item solve the word problem (are two diagrams equivalent, i.e. do they have the same interpretation) or compute a normal form for a diagram, 
\item solve the implementability problem (construct a word equivalent to a target using a library of allowed morphisms), and 
\item choose morphisms in a diagram to best approximate a more general diagram (possibly allowing the approximating diagram itself to vary). 
\end{enumerate}
Many practical questions are instances of one of these problems including computational challenges in probabilistic graphical models, quantum programming and logic \cite{coecke2011interacting,kissinger2011quantomatic}, the tensor network state approach to quantum condensed matter and quantum chemistry, parts of computational complexity theory including constraint satisfaction and counting constraint satisfaction problems, and many database operations.  

For general monoidal categories, these problems are difficult: at least $\mathsf{\#P}$-hard (1), undecidable (2), undecidable (3), and $\mathsf{NP}$-hard (4) respectively.  The news is similarly bad for approximate versions of each.  Nevertheless, given their practical significance, many tractable special cases, approximate algorithms, and heuristics exist to solve these problems in restricted cases.

Perhaps most prominent among such algorithms is the belief propagation algorithm \cite{pearl1982reverend}, and its many extensions and analogs.  From the categorical point of view, these extensions and analogs should just be the same abstract algorithm operating in different categories (e.g.\ probabilistic graphical models vs. sets and relations).  These analogies have been drawn explicitly in many areas.  For example, the connection between belief propagation and  turbo coding theory was described in \cite{mceliece1998turbo} and the connection to  survey propagation for SAT problems is explored in \cite{maneva2007new}.

We outline a general categorical form of the belief propagation algorithm for Problem 1 and look at how it specializes.  We categorified another class of algorithms for Problem 1 in previous work \cite{morton2013generalized}. 
We also touch briefly on Problems 2 and 3.  Approaches to Problem 4 often require a solution to Problem 1 as a component.

By describing algorithms in terms of monoidal categories, the common structure of problems can be understood and computational knowledge can be more readily shared across disciplines. 
Creating general tools that work for any category with suitable properties, and can be specialized automatically once the monoidal category interface of a domain was specified, would be a significant advance. 
Given the rapidly expanding universe of applied problems given categorical interpretations, such an abstraction has the potential to be as useful as convex programming or numerical linear algebra.  We only take a tiny step in this direction in the present work.

\section{
Tensor schemes and word problems} \label{sec_MCIRS}

We assume the reader is familiar with monoidal categories \cite{maclane1998categories}.  We mainly consider the strict version here.

\begin{defn}
A (finite) {\em tensor scheme} (also called a {\em monoidal signature} or {\em monoidal alphabet}) $\scrT$ is a finite set $\Ob(\scrT)$ of object variables (which must include a monoidal identity object $I$), a finite set $\Mor(\scrT)$ of morphism variables, and functions $\dom, \cod:\Mor(\scrT) \ra \Ob(\scrT)^{\ot}$.
\end{defn}

The monoidal language $\scrT^{\ot,\circ}$ generates the free monoidal category over $\scrT$. It consists of all valid morphism words that can be formed from $\Mor(\scrT)$, and identity morphisms.  Constructively, $\scrT^{\ot,\circ}$ is described as follows.
\begin{enumerate}
\item For all $A \in \Ob(\scrT)$, $\id_A$ is a word.
\item Each $f \in \Mor(\scrT)$ is a word.
\item Given words $u,u'$, $u \otimes u'$ is a word with domain $\dom(u) \otimes \dom(u')$ and codomain $\cod(u) \otimes \cod(u')$. 
\item Given words $w,w'$ with $\dom(w')=\cod(w)$, $w \circ w'$ is a word.
\end{enumerate}

Given a word in a monoidal language, we can attach an {\em interpretation} by considering it as defining a morphism in a monoidal category.  By the universal property, we can safely consider the word in the free monoidal category, which already imposes some equivalences such as $\id_A \circ f = f$ and $(f \ot g) \circ (f' \ot g') = (f \circ f') \ot (g \circ g')$.  
Two words are equivalent if they represent the same morphism in the free monoidal category.  
\begin{defn}\label{def_diagram}
An equivalence class of words in the free monoidal category over a tensor scheme is called a {\em diagram}. 
\end{defn}
Further notions of equivalence arise if we add additional relations.  So far we have no normal form for words; for example $(f \ot g) \circ (f' \ot g')$ is not preferred over the equivalent $(f \circ f') \ot (g \circ g')$. 

Because there is a coherent graphical language for the free monoidal category over a tensor scheme, the word problem is not too difficult.  
Adding adjectives (special types of monoidal categories) and relations, or fixing values, so that the category is no longer free may make it easier or harder. 
\begin{proposition}\label{prop:undec}
The word problem and implementability problem in a monoidal category over a finite tensor scheme are undecidable. 
\end{proposition}
To prove this proposition one just needs to embed a known undecidable problem in some non-free monoidal category; see \cite{morton2012undecidability} for an example of how to do this.

Proposition \ref{prop:undec} 
puts us at a distinct disadvantage as compared to computational commutative algebra, where the word problem for polynomials is always at least computable with Gr\"obner bases \cite{bachmann2007singular}. 
We call a monoidal category {\em decidable} if its word problem is decidable.  


The existence of coherent graphical languages for some types of monoidal categories means that the word problem can be reduced to graph isomorphism \cite{dixon2013open}.  Hence the word problem for the free closed category and free compact closed category over a finite tensor scheme are in $\mathsf{LOGSPACE}$ and $\mathsf{P}$ \cite{luks1982isomorphism} respectively.  Normal forms for such graphs were explored in \cite{furer1983normal,babai1983canonical}, see also \cite{mena2012trivalent}.  
One could produce normal forms for words in X-categories (for some adjective X such as ``traced''  or ``dagger'') with coherent graphical languages indirectly by this method, by transforming a word to normal form graph and then to a word again by some deterministic scheme. 

It would be preferable for some purposes to have a confluent terminating rewriting system that attached a direction to the equalities of the X-category. 
For example $(f \circ f') \ot (g \circ g') \mapsto (f \ot g) \circ (f' \ot g')$. 
Term rewriting and computing normal forms in monoidal categories is a field in its infancy; see \cite{kissinger2012pictures,mimramtowards} for some of what is known. 

A final method for the word problem is to apply a functor to a category which is complete  with respect to the category of interest but might have an easier word problem.  
Finite dimensional vector spaces over a field of characteristic zero are complete for traced symmetric monoidal categories \cite{hasegawa2008finite} and finite dimensional Hilbert spaces are complete for dagger compact closed categories \cite{selinger2011finite}. 
Thus we can potentially certify {\em inequality} of words once we have bounded dimension, for example obtaining 
numerical methods for the word problem in categories by assigning random morphisms in the category of finite-dimensional vector spaces and linear transformations.  

We assume our category is small, so the set of all morphisms from $A$ to $B$ is a {\em hom set} denoted $\Mor(A,B)$.  
We can specialize the word problem to particular hom sets such as $\Mor(I,I)$. 
In this case it is sometimes called a contraction problem (Problem 1). 
For example, in {\em semiringed} categories, the word problem for morphisms in $\Mor(I,I)$ generalizes counting constraint satisfaction problems \cite{morton2013generalized}.  

An important word problem for current purposes is determining equality of $I$-valued points, i.e.\  morphisms of type $\Mor(I,A)$ for some object $A$. The reason for this is as follows.

 We want to generalize algorithms such as belief propagation that work over the category of vector spaces and linear transformations (or a probabilistic version thereof).  
In the generalization, we can no longer assume that objects $A$ are sets with points (such as probability distributions in the classical belief propagation algorithm).  However, messages are still morphisms of type $\Mor(I,A)$ for each object $A$ and we still express the belief propagation equations in terms of equality of morphisms in each $\Mor(I,A)$.  
Deciding if two vectors are equal up to numerical tolerance becomes deciding a word problem in $\Mor(I,A)$. These messages must also be stored somehow. 

Thus for our algorithm to run efficiently, we need the word problem for $I$-valued points to be efficiently decidable and their representation to be efficient. 
Thus, we will generally assume that $I$-valued points can be stored and compared efficiently.  This can be made precise by introducing a size function. 

In the classical setting for belief propagation there is a monoid homomorphism, $\size \!:\! \Ob(\scrT)^\otimes \rae \N$,  from the free monoid generated by the objects of our tensor scheme to the natural numbers.  This sends monoidal product to multiplication (or to the addition of log dimensions if we prefer to do something closer to counting wires).  For example the dimension of a vector space fits this description.  Then for each object $A$, words in $\Mor(I,A)$ require $O(\size(A))$ storage and the word problem for $I$-valued points of $A$ is linear in the size of $A$.  


\bigskip

\subsection{Compact closed, spidered, and dungeon categories}
A compact closed category is a monoidal category with duals for objects and a compatible symmetric braiding, and string diagrams define the free compact closed category over a tensor scheme \cite{kelly1980coherence}.


\begin{defn}
A {\em spidered category} is a strict symmetric monoidal category 
 equipped with a special commutative Frobenius structure \cite{coecke2008new} $(A,m,u,\delta, \epsilon, \sigma^F)$ on each object $A$.
\end{defn}
Note that the morphisms of a spidered category need not be Frobenius, or even monoid or comonoid homomorphisms (in fact requiring this trivializes the structure). 

To such a spidered category we now add duals for objects to obtain a compact closed category with additional structure. 
We call a compact closed category which is a spidered category in a compatible way a  dungeon category.
\begin{defn}
A {\em dungeon category} is a compact closed category $(\mcC,\sigma^{\mcC},i,e)$ such that 
\begin{enumerate}[(i)]
\item Each object has a special commutative Frobenius structure $(A,m,u,\delta, \epsilon, \sigma^F)$ with $\sigma^F_{A^{(*)},A^{(*)}}=\sigma^\mcC_{A^{(*)},A^{(*)}}$, and 
\item Any two morphisms $f,g$ 
which are constructed from the identity $\id_A$, the symmetric braiding $\sigma_{A,A}$,  the Frobenius morphisms, and the dualizing cup and cap morphisms $i_A, e_A$ for $A$, and have the same domain (tensor product of zero or more or copies of $A$ and $A^*$) and the same codomain (another such tensor product), and are connected as a graph, are equal.
\end{enumerate}
\end{defn}
In other words, a ``directed spider'' morphism as in (ii) depends only on the number of inward and outward directed arrows, which way they point, and their order, with the second condition dropped up to application of $\sigma_{A,A}$.


Dungeon categories are a good setting for generalized belief propagation because we can bend wires and have spiders 
that play the role of variables in the probabilistic setting for belief propagation.

\section{Sum-product and belief propagation for contraction}\label{sec_alg_cont}
The contraction problem in semiringed categories is $\mathsf{\#P}$-hard. 
There are many cases however where the problem becomes tractable.  The main examples include trees (or diagrams where sections are merged to yield a tree) and categories for which a categorical generalization of holographic algorithms can be made to hold \cite{morton2013generalized}.


Most abstractly, the sum product algorithm \cite{kschischang2001factor} is simply the observation that when a diagram is a tree, one can perform contraction according to the tree.  
If the given diagram is not a tree, we can group nodes into a tree decomposition \cite{halin1976s} to force it to be a tree, and then run sum-product.  This is known as the {\em junction tree algorithm} \cite{lauritzen1988local} and has also been extended to the quantum case \cite{markov2008simulating}.

For actual computations, we can improve on the abstract sum-product algorithm by using an optimized message-passing version, which among other benefits permits parallelization. 

\subsection{Belief propagation in factor graphs}

First we review BP in probabilistic factor graph models with discrete variables (see e.g.\ \cite[Ch. 14]{mezard2009information}).  
The algorithm operates on a {\em factor graph}, a bipartite graph with one part discrete random variables $v \in V$ and one part factors $u \in U$.  Each factor assigns a real number to each combination of states of the variables it is connected to.  Multiplying factors and normalizing if needed gives a joint probability distribution. 

Belief propagation is a message passing algorithm.  
Each message is a probability distribution over the states one variable $v$ can take, so a vector in the associated vector space $V_v$. Each factor $f_U$ at node $u$ is a tensor in $\otimes_{v \in \nbhd(u)} V_v$. Thus a factor defines $|\nbhd(u)|$ linear maps $f_{u,v}: \otimes_{i \in \nbhd(u) \setminus v} V_i \ra V_v$, one for each $v \in \nbhd(u)$. 
First we describe how to compute messages locally at each node in the factor graph, then how to assemble these into a complete algorithm.

{\bf Messages at variables.} Suppose we have designated an edge $e$ incident on the variable as output and the rest as input.  Compute the pointwise (Hadamard) product of the incoming messages, and output it as the outgoing message along $e$.  Since we are in a probabilistic category, this Hadamard product includes a rescaling so that the message is a probability distribution.  If there are no incoming messages, output the uniform message. 

{\bf Messages at factors.} Suppose we have designated one of the edges $(u,v)$ connected to a variable $v$ incident on the factor $u$ as output and the rest as input.  Compute the tensor product of the incoming messages, apply $f_{u,v}: \otimes_{i \in \nbhd(u) \setminus v} V_i \ra V_v$, and output the result as the outgoing message along the edge to $v$.  

{\bf Resulting algorithm.} This defines a system of {\em BP equations} describing the fixed points of the update rules.  The initial messages can chosen to be uniform distributions.   When applied to a factor graph which is a tree, the algorithm converges after iterating a number of times equal to the diameter of the tree.  Choosing a root, this can be completed in two ``passes,'' leaves to root then root to leaves, updating messages only as they change.  More generally we consider a convergence threshold for the BP equations.
It is a theorem that belief propagation is exact on trees, and it can work surprisingly well even when this is not satisfied. 



\subsection{Categorical Belief Propagation}

In order for our generalized BP algorithm to operate, the category to which it is applied will need a few basic features; in particular it must be {\em reshapable}.

Suppose we have a morphism $f: \dom(f) \ra \cod(f)$ and decompositions of its domain and codomain into monoidal products of other objects: $\dom(f)=A_1 \otc A_n$ and $\cod(f) = A_{n+1} \otc A_m$. Suppose $I$ and $J$ are ordered disjoint subsets of $[m]:=\{1,2,\dots, m\}$ whose union is $[m]$. These fix  an alternate decomposition of $A_1 \otc A_m$ into two objects (with arbitrary reordering), $D= \ot_{i \in I} A_i$, $C=\ot_{j \in J}A_j$.  We require that there exists a unique morphism $r_{I,J}(f):D \ra C$ called the {\em reshaping of $f$ from $I$ to $J$} determined from this data.  The reshaping does not depend on the route, that is for any $I',J'$  we have $r_{I,J}r_{I'J'}=r_{I,J}$.  If this works for all morphisms, we say the category is {\em reshapable}.

We now assume that we are working in a dungeon category (which is therefore reshapable) in which $I$-valued points can be stored and compared efficiently (e.g.\ space and time complexity linear in some ``size'' monoid homomorphism).   
This will guarantee the efficiency of belief propagation by an argument analogous to the classical case.

\begin{comment}
\begin{example}[reshaping in $\Mat_\FF$: transpose, vectorization, matricization]

\end{example}

\begin{example}[multilinear algebra: FinVect]

\end{example}

\begin{example}[FinRel and $\N$FinRel]

\end{example}

\begin{example}[Spivak's Sch]

\end{example}

\begin{example}[Spiders]

\end{example}
\end{comment}

\subsubsection{Bipartite with one part spiders}
Assume we have a dungeon category (a strict compact closed category 
 equipped with a compatible special commutative Frobenius structure on each object). 
This allows for 
the special spider morphisms to generalize and axiomatize the role played by variables in the traditional belief propagation algorithm for factor graphs.  
Consider a bipartite diagram with one part consisting of spiders (generalized variables) and the other arbitrary morphisms (generalized factors).   
Both spiders and factors are reshaped before composition.  The reshaping of the spiders simply leads to a Hadamard (elementwise) product in the case of vector spaces and linear transformations. 

The BP algorithm is expressed in terms of messages.  When objects are sets, messages are elements of the set.  For example, if objects are vector spaces, the messages will be vectors.
In general messages for object $A$ are $I$-valued points (elements of $\Mor(I,A)$). 
For example any probability distribution can be expressed this way as a stochastic matrix applied to the unit Frobenius morphism (the unit  ``creates'' a variable with a uniform distribution).


Using the compact closed structure to reshape morphisms, we still have that each ``factor'' morphism at node $u$ defines $|\nbhd(u)|$ different morphisms $f: \otimes_{i \in \nbhd(u) \setminus v} V_i \ra V_v$, one for each $v \in \nbhd(u)$. 
First we describe how to compute messages locally at each node, then how to assemble these into a complete algorithm.

{\bf Messages at spiders.} 
Apply the reshaped spider to incoming messages, and output the result as the outgoing message.  If there are no incoming messages, treat the spider as a Frobenius unit. 

{\bf Messages at ``factor'' morphisms.} Suppose we have designated one of the wires incident on the factor $f$ as output and the rest as input.  Compute the monoidal product of the incoming messages, apply the reshaped $f$, and output the result as the outgoing message. 

The system of  BP equations are now equalities of $I$-valued morphisms, describing the fixed points of the update rules.  The initial messages can be chosen to be units at the spiders.  The nice behavior of the algorithm on diagrams which are trees should be preserved, once suitable definitions are made so that we can talk about convergence.






A  spider is just a special kind of morphism.  To get the {\bf general bipartite} version, replace the message procedure at spiders with another copy of the factor message procedure.

\subsection{Examples}\label{sec_examples}




In the category $\BoolRel$ each object is a two element set or a monoidal (here Cartesian) product of such. Morphisms are relations.  Categorical belief propagation becomes the survey propagation algorithm for solving constraint satisfaction problems.

Augmenting the category $\FinRel$ of finite sets and relations with a positive integer multiplying each element of a relation yields the semiringed category $\NFinRel$ (discussed in the context of counting constraint satisfaction problems such as computing partition functions in \cite{morton2013generalized}).  This category also corresponds approximately to sufficient statistics in the analysis of contingency tables \cite{agresti2013categorical}.

We can add additional flexibility to $\FinRel$ to obtain database categories \cite{spivak2012functorial}.  Then categorical belief propagation becomes a query planning algorithm.  







The category in which numerical linear algebra takes place is vector spaces and linear transformations, where the vector spaces are augmented by orthonormal bases which define spiders. 
Composition is matrix multiplication, and tensor product is Kronecker product.  
Consider the {\em dual numbers} over a field $\FF$,  $\mathbb{D} = \FF[\epsilon]/\langle \epsilon^2 \rangle$.  Reverse-mode automatic differentiation  is categorical belief propagation in the category of vectors and matrices over $\mathbb{D}$.

We also conjecture \cite{critch2012polynomial} that algorithms commonly used in quantum condensed matter physics such as DMRG \cite{PhysRevLett.69.2863} and its many extensions can be considered as an instance of the categorical belief propagation algorithm.













\bibliographystyle{eptcs}
\bibliography{references}
\end{document}